\documentclass[12pt]{article}
\usepackage{amssymb}
\makeatletter
\date{}

\newtheorem{theorem}{Theorem}[section]

\newtheorem{proposition}[theorem]{Proposition}

\newtheorem{problem}[theorem]{Problem}

\newcommand{\edim}{{\rm e}$-${\dim}}

\newcommand{\z}{{\Bbb Z}}

\newcommand{\lo}{\longrightarrow}

\newcommand{\black}{{\blacksquare}}

\begin{document}

\title{A dimensional property of Cartesian product}

\author{  Michael Levin\footnote{the author was supported by ISF grant 836/08}}

\maketitle
\begin{abstract} 
We show that the Cartesian product of  three hereditarily infinite dimensional
compact metric spaces is never hereditarily infinite dimensional.
It is quite surprising that the proof of this fact (and this is the only
proof known to the author) essentially relies on algebraic topology.
\bigskip
\\
{\bf Math. Subj. Class.:}  55M10 (54F45 55N45)
\end{abstract}
\begin{section}{Introduction}\label{intro}
Throughout this paper we assume that  maps are continuous and
spaces are separable metrizable. We recall that a  compactum
means a compact metric space. By the dimension $\dim X$  of a space $X$ we
assume the covering dimension.

An infinite dimensional  compactum $X$ is said to be hereditarily infinite 
dimensional 
if every (non-empty)  closed subset of $X$ is either $0$-dimensional 
or infinite dimensional. Hereditarily infinite dimensional compacta were
first constructed by Henderson \cite{henderson}, for 
 related results and simplified constructions
 see \cite{rubin}, \cite{walsh}, \cite{pol}, \cite{levin1}
\cite{levin2}. The main result of this paper is: 
\begin{theorem}
\label{main-th}
Let $n>0$ be an integer and $X_i$, $1\leq i \leq n+2$, hereditarily infinite dimensional compacta.
Then the product $Z=\Pi_{i=1}^{n+2} X_i$ contains an $n$-dimensional closed subset.
In particular, the product of three hereditarily infinite dimensional compacta
is never hereditarily infinite dimensional.
\end{theorem}

Let us  note that in general $Z$ in Theorem \ref{main-th} does not contain 
finite dimensional subspaces of arbitrarily large dimension. Indeed, consider
the Dydak-Walsh compactum  $X$ \cite{dydak-walsh} having the following properties:
 $\dim X =\infty$,  $\dim_\z X =2$ and 
$\dim_\z X^n =n+1$ for every positive integer $n$. 

We remind
that for an abelian group $G$ the cohomological dimension $\dim_G X $ of a space $X$
is the smallest integer $n$ such that  the Cech cohomology $H^{n+1} (X, A; G)$ 
vanish for every closed subset $A$ of $X$.
Clearly $\dim_G X \leq \dim X$ for every abelian group $G$.
By the classical result of  Alexandroff 
$\dim X =\dim_\z X$ if $X$ is finite dimensional. Alexandroff's result was
extended by  Ancel \cite{ancel}
who showed that 
$\dim X= \dim_\z X$ if $X$ is a compact C-space. 
We remind that a space $X$ is a C-space if
for any infinite sequence ${\cal U}_i$  of open covers $X$ there is an open  cover
${\cal V}$ of $X$ such that
  ${\cal V}$ splits into the union ${\cal V}= \cup_i {\cal V}_i$
of families ${\cal V}_i$ of disjoint  sets   such that
${\cal V}_i$ refines ${\cal U}_i$. 

Thus the Dydak-Walsh compactum $X$
is not a C-space. R. Pol \cite{pol-c-space} (see also \cite{levin2})
 showed that a compactum which is not a C-space contains a hereditarily 
 infinite dimensional closed subset. Hence, replacing  $X$
 by its  hereditarily infinite dimensional closed subset, we may assume that
 $X$ is hereditarily infinite dimensional. Since $\dim_\z X^{n+2}=n+3$,
 we deduce from  Alexandroff's theorem that $X^{n+2}$ does not contain
 finite dimensional subsets of $\dim> n+3$. Moreover, $X^{n+2}$ does not
 contain compact subsets of $\dim=n+3$. Indeed, if $F$ is a finite dimensional
  closed subset of  $X^{n+2}$ then, since $X$ is hereditarily infinite dimensional,
  the projection of $p : F \lo X^{n+1}$ is $0$-dimensional. By a result of Dranishnikov
  and Uspenskij
  \cite{d-u} a 0-dimensional map  of compacta
  cannot lower cohomological dimensions and hence
  $\dim F =\dim_\z  F \leq \dim_\z X^{n+1} =n+2$.  
  
  This example together with Theorem \ref{main-th} suggest

  \begin{problem}
  Does the compactum $Z$ in Theorem \ref{main-th} always contain
  a closed subset of $\dim=n+1$? of $\dim=n+2$? a subset of $\dim=n+1$?
  of $\dim =n+2$? of $\dim=n+3$ ?
  \end{problem}
  
  Note that  Theorem \ref{main-th} implies that there are two hereditarily infinite 
  dimensional compacta whose product is not hereditarily infinite dimensional.
  Indeed, let $X_1$, $X_2$ and $ X_3$ be hereditarily infinite dimensional compacta. 
   If $X_1\times X_2$ is hereditarily infinite dimensional
  then, by Theorem \ref{main-th}, $(X_1 \times X_2) \times X_3$ is not hereditarily infinite
  dimensional. This observation suggests
  
  \begin{problem}
  Do there exist two hereditarily infinite dimensional compacta
  whose product is also hereditarily infinite dimensional? Does there
  exist a hereditarily infinite dimensional compactum whose
  square is hereditarily infinite dimensional?
  \end{problem}
  
  It is quite surprising that the proof of Theorem \ref{main-th}
  essentially relies on algebraic topology.
  It would be interesting to find an elementary direct proof of Theorem \ref{main-th}.
 \end{section}

 \begin{section}{ Proof of Theorem \ref{main-th}}
 Let us recall basic  definitions and results in Extension Theory and Cohomological
 Dimension that will be used in the proof.
 
 The extension dimension of a space $X$ is said to be  dominated by 
 a CW-complex $K$, written $\edim X \leq K$, if every map $f : A \lo K$ from
 a closed subset $A$ of $X$ extends over $X$. We also write $\edim X > K$
 if the property $\edim X \leq K$ does not hold.  Note
 that the property $\edim X \leq K$ depends only on the homotopy type
 of $K$.  
 The covering and cohomological dimensions can be 
 characterized by the following extension properties: 
 $\dim X \leq n$ if and only if the extension dimension of $X$ is  dominated by
 the $n$-dimensional sphere $S^n$ and 
 $\dim_G X \leq n$ if and only if the extension dimension of 
 $X$ is  dominated by the Eilenberg-Mac Lane complex
 $K(G,n)$. The extension dimension shares many properties of covering dimension.
 For example: 
 if $\edim X \leq K$ then for every $A \subset X$ we have
 $\edim A \leq K$, and if $X$ is a countable union of closed subsets whose
 extension dimension is dominated by $K$ then $\edim X \leq K$. In the proof of 
 Theorem  \ref{main-th} we will also use the following facts.

 \begin{theorem}
 \label{olszewski}
 {\rm \cite{olszewski}}
 Let $K$ be a countable CW-complex and $A$ a subspace of a compactum $X$
 such that $\edim A \leq K$. Then there is a $G_\delta$-set $ A' \subset X$
 such that $A \subset A'$ and $\edim A' \leq K$.
 \end{theorem}
 
 \begin{theorem}
 \label{dranishnikov-splitting}
{\rm  \cite{dranishnikov-splitting}}
 Let $K$ and $L$ be countable CW-complexes and  $X$ a compactum such that
 $\edim X \leq K*L$. Then $X$ decomposes into subspaces $X=A \cup B$ 
 such that $\edim A  \leq K$ and $\edim B \leq L$.
 \end{theorem}
 
 \begin{theorem}
 \label{levin}
 {\rm \cite{levin}}
 Let $f : X \lo Y$ be a map of compacta and let
 $K$ and $L$ be  countable CW-complexes such that 
 $\edim Y \leq K$ and $\edim f^{-1}(y) \leq L$ for every
 $y \in Y$. Then $\edim X \leq K*L$. In particular, if for a compactum $Z$
 we have $\edim Z \leq L$ then $\edim Y \times Z \leq K *L$.
 \end{theorem}
 
 \begin{theorem}
 \label{dranishnikov-characterization}
 {\rm \cite{dranishnikov-characterization}}
 Let for a compactum $X$ and a CW-complex $K$ we have $\edim X \leq K$.
 Then $\dim_{H_n(K)} X \leq n$ for every $n\geq 1$.
 \end{theorem}
  
 By $\z_p$ we denote the $p$-cyclic group and by $\z_{p^\infty}=dirlim \z _{p^k}$ 
 the $p$-adic circle.
 
 \begin{theorem}
 \label{kuzminov-bockstein}
{\rm  \cite{kuzminov-bockstein}}
 Let $p$ be a prime and $X$ and $Y$  compacta.
 Then $\dim_{\z_p} X\times Y =\dim_{\z_p} X + \dim_{\z_p} Y$.
 \end{theorem}
 
 \begin{theorem}
 \label{dydak-sum}
 {\rm \cite{dydak-union}}
Let $p$ be a prime and $X=A \cup B$ a decomposition  of a compactum $X$.
Then $\dim_{\z_p} X \leq \dim_{\z_p} A + \dim_{\z_p} B  +1$.
\end{theorem}

  For an abelian group $G$ 
 we always  assume that a Moore space  $M(G, n)$ of type $(G,n)$ is a CW-complex and 
 $M(G,n)$ is simply connected  if $n>1$.
 Note that $M(G,n)$ is defined uniquely (up to homotopy equivalence) for $n>1$
 \cite{hatcher}. 
 Recall that for  CW-complexes $K$ and $L$ the join $K*L$ is
  homotopy equivalent to the suspension $\Sigma (K \wedge L)=S^0 *(K \wedge L)$ and
 $K * L$  is simply
 connected if at least one of the complexes $K$ and $L$ is connected.
 Then it follows from the K${\rm \ddot{ u}}$nneth formula that for 
 distinct primes $p$ and $q$:\\
 
 (i)  $M(\z_p, 1) * M(\z_q, 1)$ is contractible;
 
 (ii)  $M(\z_{p^\infty}, 1) * M(\z_q, 1)$ is contractible;
 
 (iii)  $M(\z_{p^\infty}, 1) * M(\z_{q^\infty}, 1) $ is contractible;
 
 (iv) $\Sigma^n M(\z_p,1)=S^{n-1}*M(\z_p, 1)$ is a Moore space $M(\z_p, n+1)$;
 
(v)  $M(\z_{p^\infty}, 1) * M(\z_p, n)$ is a Moore space $M(\z_p, n+3)$.
\\\\
We  say that a compactum $X$ is {\em reducible} at a prime $p$ if
there is  a non-zero-dimensional closed subset $F$ of $X$ with
$\edim F \leq M(\z_p, 1)$  and $\edim F \leq M(\z_{p^\infty},1)$, and
 we say that $X$ is {\em unreducible} at  $p$  otherwise.

\begin{proposition}
\label{singular-value}
Let $X$ be a hereditarily infinite dimensional compactum. Then 
$X$ is unreducible at  at most one prime.
\end{proposition}
{\bf Proof.} Aiming at a contradiction 
assume $X$ is unreducible at two  distinct primes $p$ and $q$.
By (i) we have $\edim X \leq M(\z_p,1)*M(\z_q,1)$ and hence, by  
Theorem \ref{dranishnikov-splitting}, the compactum $X$ decomposes into
$X=A \cup B$ with $\edim A \leq M(\z_p,1)$ and $\edim B \leq M(\z_q,1)$
and, by Theorem \ref{olszewski}, we may assume that $B$ is $G_\delta$ and
$A$ is $\sigma$-compact.  

If $\dim A >0$ then $A$ contains  a non-zero-dimensional compactum 
 $F \subset A$  and clearly $F$ is hereditarily infinite dimensional 
and $\edim F \leq M(\z_p, 1)$. 

If $\dim A \leq 0$ then replacing $A$ by a bigger  
 $0$-dimensional  $G_\delta$-subset  of $X$
 we may assume that 
$B$ is $\sigma$-compact. Since $X$ is 
infinite dimensional we have that  $\dim B>0$ and  
hence $B$ contains a non-zero-dimensional compactum
 $F \subset B$. Clearly $F$ is hereditarily infinite dimensional 
and $\edim F \leq M(\z_q, 1)$. 

Thus without loss of generality we may assume that $X$ contains
a hereditarily infinite dimensional compactum $F$ with
$\edim F \leq M(\z_p,1)$. By (iii) we have that 
$\edim F \leq M(\z_{p^\infty},1)*M(\z_{q^\infty},1)$.  Then  using 
the above reasoning  we can replace  $F$ by  a hereditarily infinite dimensional 
closed subset of $F$ and assume, in addition, that
the extension dimension of $F$ is dominated by at least one the complexes
$M(\z_{p^\infty},1)$ or $M(\z_{q^\infty},1)$. 

If $\edim F \leq M(\z_{p^\infty},1)$ then $X$ is reducible at  $p$  
and we are done. If $\edim F \leq M(\z_{q^\infty},1)$ then, by (ii), 
we have that $\edim F \leq M(\z_{p^\infty},1) * M(\z_q,1)$ and once again
by the reasoning described above one can replace $F$ by its
closed hereditarily infinite dimensional subset with 
the extension dimension dominated by at least one of the complexes
$ M(\z_{p^\infty},1)$  or $ M(\z_q,1)$. This implies that $X$ is reducible  at at least one of 
the primes $p$ and $q$  and the proposition follows. 
 $\black$
\\

  Now we are ready to prove Theorem \ref{main-th}.
 \\
 {\bf Proof of Theorem \ref{main-th}.} By Proposition \ref{singular-value} there is a prime $p$
 such that  every $X_i$ is reducible at $p$. Hence for every $i$ there is
 a hereditarily infinite dimensional compactum $F_i \subset X_i$ such that
 $\edim F_i \leq M(\z_p,1)$ and $\edim \leq M(\z_{p^\infty},1)$. Then, 
 by Theorem \ref{dranishnikov-characterization}, $\dim_{\z_p} F_i=1$ and, 
 by Theorem \ref{kuzminov-bockstein}, $\dim_{\z_p} F =n+2$ for $F=F_1 \times \dots \times F_{n+2}$.
 
 On the other hand, by Theorem \ref{levin}, 
 $\edim F \leq K=M(\z_{p^\infty},1)*\dots* M(\z_{p^\infty},1)*M(\z_p,1)$
 (the join of $M(\z_p,1)$ and  $n+1$ copies of $ M(\z_{p^\infty},1)$). By (v) and (iv)
 we have that $K=M(\z_p, 3n+4)= S^{3n+2}*M(\z_p,1)$. 
 Then, 
 by Theorem \ref{dranishnikov-splitting}, $F$ splits into $F=A \cup B$
 such that  $\edim A \leq M(\z_p, 1)$ and $B$ is finite-dimensional. In addition,
 we may assume by Theorem \ref{olszewski} that $B$ is $G_\delta$ and  
 $A$ is $\sigma$-compact. Then, by Theorem \ref{dranishnikov-characterization},
 the property $\edim A \leq M(\z_p,1)$ implies $\dim_{\z_p} A \leq 1$. 
  Again
 by Theorem \ref{olszewski}, we can replace $A$ by a bigger $G_\delta$ subset 
 of $F$ and assume  that  $\dim_{\z_p} A  \leq 1$ and $B$ is 
  finite-dimensional and $\sigma$-compact.
  
Then,
 by Theorem \ref{dydak-sum},  we have 
 $n+2=\dim_{\z_p} F \leq \dim_{\z_p} A + \dim_{\z_p} B +1\leq \dim_{\z_p} B +2$ and
 hence  $\dim_{\z_p} B \geq n$. Thus $\dim B\geq n$ and, since $B$ is 
 finite dimensional and $\sigma$-compact,  $B$ contains an $n$-dimensional
 compact subset. The theorem is proved. $\black$

\end{section}

Michael Levin\\
Department of Mathematics\\
Ben Gurion University of the Negev\\
P.O.B. 653\\
Be'er Sheva 84105, ISRAEL  \\
 mlevine@math.bgu.ac.il\\\\

\begin{thebibliography}{99}
\bibitem{ancel}
Ancel, Fredric D. The role of countable dimensionality in the theory of cell-like relations.  
Trans. Amer. Math. Soc.  287  (1985),  no. 1, 1–-40. 

\bibitem{dranishnikov-characterization}
Dranishnikov, A. N. Extension  of mappings into CW-complexes. (Russian)  
Mat. Sb.  182  (1991),  no. 9, 1300--1310;  translation in  Math. USSR-Sb.  74  (1993),  
no. 1, 47–-56.

\bibitem{dranishnikov-splitting}
Dranishnikov, A. N. On the mapping intersection problem.  
Pacific J. Math.  173  (1996),  no. 2, 403–-412.

\bibitem{d-u}
 Dranishnikov, A.N.; Uspenskij, V.V. Light maps and extensional dimension, 
Topology Appl. 80
(1997) 91–-99.

\bibitem{dydak-union}
Dydak, Jerzy Cohomological dimension and metrizable spaces. II.  
Trans. Amer. Math. Soc.  348  (1996),  no. 4, 1647–-1661.

\bibitem{dydak-walsh}
Dydak, Jerzy; Walsh, John J.
 Infinite-dimensional compacta having 
cohomological dimension two: An
application of the Sullivan conjecture, Topology 32 (1) (1993) 93-–104.

\bibitem{hatcher}
Hatcher, Allen Algebraic topology. Cambridge University Press, Cambridge, 2002. xii+544 pp. 
ISBN: 0-521-79160-X.

\bibitem{henderson}

Henderson, David W. An infinite-dimensional compactum with 
no positive-dimensional compact subsets—a simpler construction.  
Amer. J. Math.  89  1967 105–-121. 

\bibitem{kuzminov-bockstein}
 Kuzminov, V.I. Homological dimension theory, Russian Math. Surveys 23, issue 5 (1968), 1--45.

\bibitem{levin1}
Levin, Michael A short construction of hereditarily infinite-dimensional compacta.  
Topology Appl.  65  (1995),  no. 1, 97–-99.

\bibitem{levin2}
Levin, Michael Inessentiality with respect to subspaces.  Fund. Math.  
147  (1995),  no. 1, 93–-98.

\bibitem{levin}
Levin, Michael On extensional dimension of maps.  
Topology Appl.  103  (2000),  no. 1, 33–-35.

\bibitem{pol}
Pol, Roman
Selected topics related to countable-dimensional metrizable spaces.  
General topology and its relations to modern analysis and algebra, VI (Prague, 1986),  421–-436,
Res. Exp. Math., 16, Heldermann, Berlin, 1988. 

 \bibitem{pol-c-space}
 Pol, Roman On light  mappings without perfect fibers on compacta.  
 Tsukuba J. Math.  20  (1996),  no. 1, 11–-19. 
 
 \bibitem{olszewski}
 Olszewski, Wojciech
  Completion theorem for cohomological dimensions, Proc. Amer. Math. Soc. 123
(7) (1995) 2261–-2264.
\bibitem{rubin}
Rubin, Leonard R.
Hereditarily strongly infinite-dimensional spaces.
Michigan Math. J. 27 (1980), no. 1, 65–-73. 

\bibitem{walsh}
Walsh, John J. Infinite-dimensional  compacta containing 
no n-dimensional ($n\geq 1$) subsets.  
Topology  18  (1979), no. 1, 91–-95.

\end{thebibliography}
\end{document}